\relax 
\global \edef \rmkrrred {{1.2}}
\global \edef \rmkdepth {{1.6}}
\global \edef \exnotCM {{1.7}}
\global \edef \exaci {{1.8}}
\global \edef \propmuplusone {{1.9}}
\global \edef \propxyl {{1.11}}
\global \edef \propxyaal {{1.12}}
\global \edef \propdepthone {{2.2}}
\global \edef \exjarrah {{3.3}}
\global \edef \thmr {{4.1}}
\global \edef \reductionno {{4.2}}
\global \edef \reduction {{4.5}}

\magnification 1100
\overfullrule=0pt
\parskip=3pt
\baselineskip=15.0pt

\font\bigfont=cmr12

\font\bbigfont=cmbx12

\font\goth=eufm10

\def\qed{\hfill\hbox{\vrule\vtop{\vbox{\hrule\kern1.7ex\hbox{\kern1ex}
  \hrule}}\vrule}\vskip1.5ex}
\def\eqed{\hbox{\quad
  \hbox{\vrule\vtop{\vbox{\hrule\kern1.7ex\hbox{\kern1ex}\hrule}}\vrule}}}
\def \r#1{{\hbox{$\widetilde{#1}$}}}

\def\lt{\hbox{lt}\;}
\def\Ass{\hbox{Ass}\;}

\def\mm{\hbox{\goth m}}

\def\rrelt{\underline X ^{\underline\mu}}
\def\rselt{\underline X ^{\underline\mu'}}
\def\gelt{\underline X^{\underline\nu}}

\def\overset#1#2{{\mathop{\buildrel {#1} \over {#2}}}}

\font\tenmsb=msbm10
\font\sevenmsb=msbm7
\font\fivemsb=msbm5
\newfam\msbfam
\textfont\msbfam=\tenmsb
\scriptfont\msbfam=\sevenmsb
\scriptscriptfont\msbfam=\fivemsb

\def\hexnumber#1{\ifcase#1 0\or1\or2\or3\or4\or5\or6\or7\or8\or9\or
	A\or B\or C\or D\or E\or F\fi}

\mathchardef\subsetneq="2\hexnumber\msbfam28

\newcount\sectno \sectno=0
\newcount\thmno \thmno=0
\def \section#1{\bigskip\bigskip
	\global\advance\sectno by 1 \global\thmno=0
	\noindent{\bf \the\sectno. #1}}
\def \thmline#1{\vskip 6pt
	\global\advance\thmno by 1
	\noindent{\bf #1\ \the\sectno.\the\thmno:}\ \ %
	\bgroup \advance\baselineskip by -1pt \it
	\abovedisplayskip =4pt
	\belowdisplayskip =3pt
	\parskip=0pt
	}
\def \dthmline#1{\vskip 6pt
	\global\advance\thmno by 1
	\noindent{\bf #1\ \the\sectno.\the\thmno:}\ \ }

\def \endb{\egroup}

\def \remark{\dthmline{Remark}}
\def \question{\thmline{Question}}
\def \example{\dthmline{Example}}
\def \examples{\dthmline{Examples}}

\def \prop{\thmline{Proposition}}
\def \proof{\smallskip\noindent {\sl Proof:\ \ }}



\def\label#1{\unskip\immediate\write\isauxout{\noexpand\global\noexpand\edef\noexpand#1{{\the\sectno.\the\thmno}}} 
{\global\edef#1{\the\sectno.\the\thmno}}%
\unskip}

\newwrite\isauxout
\openin1\jobname.aux
\ifeof1\message{No file \jobname.aux}
          \else\closein1\relax\input\jobname.aux
          \fi
\immediate\openout\isauxout=\jobname.aux
\immediate\write\isauxout{\relax}

\def\today{\ifcase\month\or January\or February\or March\or
     April\or May\or June\or July\or August\or September\or
     October\or November\or December\fi
     \space\number\day, \number\year}

\ %
\font\footfont=cmr8
\font\footitfont=cmti8
\smallskip
\centerline{\bbigfont Notes on the behavior of the Ratliff-Rush filtration}
\bigskip
\centerline{\bigfont Maria Evelina Rossi and Irena Swanson}
\unskip\footnote{ }{{\footitfont 2000 Mathematics Subject Classification.}
13A30, 13B22, 13P10
}
\bgroup
\baselineskip=5pt 
\unskip\footnote{ }{{\footitfont Key words and phrases.}
\footfont Ratliff-Rush ideals, superficial elements,
Gr\"obner bases, initial ideals, associated primes,
number of generators,
(Ratliff-Rush) reduction number.}

\egroup
\bigskip

\bgroup
\narrower\narrower
\baselineskip=10pt
\noindent
{\bf Abstract.}
We establish new classes of Ratliff-Rush closed ideals
and some pathological behavior of the Ratliff-Rush closure.
In particular,
Ratliff-Rush closure does not behave well under passage modulo 
superficial elements,
taking powers of ideals,
associated primes,
leading term ideals,
and the minimal number of generators.
In contrast,
the reduction number of the Ratliff-Rush filtration behaves better:
it preserves some information on the reduction number of the ideal.

\egroup
\bigskip

Let $I$ be an ideal in a Noetherian ring $R$.
  From the maximal
condition it follows that there exist ideals $\r{I}$ in $R$ maximal
with respect to the condition
$$
\r{I}^n = I^n \hbox{\ for all large\ } n.
$$
Ratliff and Rush proved in [RR, Theorem 2.1],
that if $I$ is a regular ideal (i.e., it contains a non-zerodivisor),
then there exists a unique largest such $\r I$,
which can be presented in terms of $I$ as follows:
$$\r{I}:= \bigcup_{n\ge1} (I^n : I^{n-1}).$$

When $I$ is not regular,
this fails.
For example,
let $F$ be a field, $X$ an indeterminate over $F$,
and $R = F[X]/(X^2)$.
Then the ideal $I = XR$ is not regular,
the (unique) largest ideal $\r I$ for which
$\r{I}^n = I^n$ for all large $n$ is $\r I = XR$,
but
$\bigcup_{n\ge1} (I^n : I^{n-1}) = R$.
In the sequel we always assume that $I$ contains a non-zerodivisor.

The ideal $\r I$ is called the {\it Ratliff-Rush ideal associated with}
$I$ (see [HJLS]) or the {\it Ratliff-Rush closure of} $I$.
A regular ideal $I$
for which $\r I =I$ is called {\it Ratliff-Rush closed}.

In this paper we analyze various properties that
Ratliff-Rush closures of ideals do or do not satisfy.
First of all,
clearly $\r{\r I} = \r I$.
Also,
if $a$ is a non-zerodivisor such that for an ideal $I$,
$I : a = I$,
then $\r{aI} = a \r I$.
However,
it is well-known that the Ratliff-Rush closure is not a ``closure''
in the usual sense.
Namely,
if $I \subseteq J$,
it need not follow that $\r I \subseteq \r J$.
An example, taken from [HJLS], is
$I=(X^3, Y^3)$ and $J=(X^4,X^3Y,XY^3,Y^4)$ in the polynomial ring
$R=F[X,Y]$ over a field $F$.
Then $X^2Y^2 \in \r J \setminus \r I$.

Nevertheless,
the Ratliff-Rush closure of ideals is a good operation
with respect to many properties,
it carries information about associated primes of powers of ideals,
about zerodivisors in the associated graded ring,
preserves the Hilbert function of zero-dimensional ideals, etc.
In the different sections we review some of the known properties,
and examine the behavior of the Ratliff-Rush closure
with respect to some other properties.
Section~1 examines the behavior of the Ratliff-Rush closure
on powers of a fixed ideal.
We present many examples illustrating exceptional behavior.
We also establish new classes of ideals
for which all the powers are Ratliff-Rush closed.
In Section~2 we show that the Ratliff-Rush closure
has no natural corresponding notion of ``superficial elements",
so that in general one cannot drop dimension in arguments
using the Ratliff-Rush closure.
In Section~3 we show that Ratliff-Rush closure also does not
behave well with respect to Gr\"obner bases, leading term ideals,
taking the sets of associated primes,
and the minimal number of generators.
This answers a questions of Heinzer et al. from [HJLS].
Section~4 shows some positive results for the Ratliff-Rush filtration,
namely that the Ratliff-Rush reduction number behaves well.

It is easy to see that $I \subseteq \r I$ and that an element of
$(I^n : I^{n-1})$ is integral over $I$.
Hence for all regular ideals
$I$,
$$I \subseteq \r{I} \subseteq \overline{I} \subseteq \sqrt{I}.
$$
Thus there exist many ideals which are Ratliff-Rush ideals,
for example, all radical and all integrally closed regular ideals.

The main thrust of our analysis is the comparison
of the behaviors of the Ratliff-Rush and integral closures of ideals.
We illustrate the different types of behaviors with examples.
Most of our examples are monomial ideals in polynomial rings over fields.
By definition the Ratliff-Rush closure of a monomial ideals is a 
monomial ideal,
and this makes some computations easier.
While there exist algorithms for computing the integral closures of ideals,
there exist no such algorithms for the Ratliff-Rush closures.
Namely,
to compute $\cup_n (I^{n+1} : I^n)$,
of course there exists a positive integer $N$
such that
$\cup_n (I^{n+1} : I^n) = I^{N+1} : I^N$.
However,
it is not clear how big this $N$ is.
Just because
$I^{n+1} : I^n = I^{n+2} : I^{n+1}$,
it does not imply that $I^{n+1} : I^n = I^{n+3} : I^{n+2}$
(see Example~\exaci \ below).
Thus our computations of the Ratliff-Rush closures,
or verifications that an ideal is Ratliff-Rush closed,
are in general laborious.
In a few cases we establish the Ratliff-Rush closure property indirectly:
we establish that the associated graded ring of the ideal has positive depth,
in which case all the powers of the ideal are Ratliff-Rush closed
by [RR, (2.3.1)], or we prove that the ideal is integrally closed.
But here is an ideal which is Ratliff-Rush closed
but not integrally closed:
let $I$ be the ideal $(X^4, X^3Y, X^2Y^2, Y^4)$
in the polynomial ring $R = F[X,Y]$,
with $F$ a field.
Then $I$ is not integrally closed as $XY^3$ is integral over $I$
but not in $I$.
However,
$I$ is Ratliff-Rush closed.
Namely,
by degree count every element of the Ratliff-Rush closure has
to have total degree at least $4$.
As $(X,Y)^5 \subseteq I$,
it suffices to verify that $XY^3$ is not in the Ratliff-Rush closure of $I$.
But $XY^3 (Y^4)^n$ is not an element of $I^{n+1}$ for any positive integer $n$,
which proves the claim.

On the other hand,
the Ratliff-Rush closure of $(X^4, X^3Y, XY^3, Y^4)$ is $(X,Y)^4$.
These two examples show that the Ratliff-Rush closure of monomial ideals
need not exhibit the integral closure's property
of the convexity of the corresponding Newton polytopes.

\vfill\eject
\section{Ratliff-Rush closure and powers of an ideal}

In this section we prove several properties of the Ratliff-Rush closures
of powers of ideals.
For example,
all the high powers are Ratliff-Rush closed,
but many pathologies occur for low powers.
We present examples of these pathologies.
The last half of the section is about ideals
all of whose powers are Ratliff-Rush closed.

Ratliff and Rush proved that for a regular ideal $I$,
for all large $n$, $(\r I)^n = I^n$.
In particular when $I$ is a zero-dimensional ideal,
$\r I$ is the largest ideal which has the same Hilbert polynomial as~$I$.
Furthermore,
high powers of an arbitrary regular ideal are Ratliff-Rush closed:

\remark
Let $I$ be a regular ideal. Then for all large $n$,
$I^n = \r {I^n}$.
(See [RR].)
%

A somewhat simpler computation of the Ratliff-Rush closures is achieved as follows:

\remark
\label{\rmkrrred}
$\r {I^n} = \bigcup_k (I^{n+k} : I^k)$.
Furthermore,
if $J=(a_1,\dots, a_d)$ is a reduction of $I$,
then $\r {I^n} = \cup_k (I^{n+k} : (a_1^k,\dots, a_d^k))$.
(See [RR].)

%
%

This all shows that the Ratliff-Rush closures of high powers of 
ideals behave well.
In particular,
for all (regular) ideals, sufficiently large powers are Ratliff-Rush closed.
(If we admit non-regular ideals,
any non-radical nilpotent ideal is a counterexample of this property.)
However,
the Ratliff-Rush closures of powers of (regular) ideals
do exhibit some pathologies.
We show examples below of these pathologies,
namely of Ratliff-Rush closed ideals whose powers are not Ratliff-Rush closed,
and of ideals which are not Ratliff-Rush closed but whose powers are.
In the last half of this section we analyze some classes of ideals 
all of whose powers
are Ratliff-Rush closed.

\example
There are integrally closed ideals (even prime ideals),
and so Ratliff-Rush closed,
whose second power is not Ratliff-Rush closed,
even if $I$ is a maximal ideal of a local Cohen-Macaulay ring.
Here is a simplification of an example due to H.\ J.\ Wang:
let $I$ be the maximal ideal of the Cohen-Macaulay local ring
$R=F[[X,Z,U]]/(Z^2,ZU,XZ-U^3)$,
where $F$ is a field and $X, Z$ and $U$ indeterminates over $F$.
As $X$ is a non-zerodivisor in $R$,
$I$ is a regular ideal.
In this case $I$ is Ratliff-Rush closed but $I^2$ is not.
In fact, $I$ is a prime ideal 
and so Ratliff-Rush closed,
but $I^3:I$ contains $Z$ which is not in $I^2$, and hence $I^2 \neq \r{I^2}$.
(Wang's original example was
$R=F[[X,Y,Z,U,V]]/(Z^2,ZU,ZV,UV,XZ-U^3,ZY-V^3)$.)

There are even monomial Ratliff-Rush closed ideals in polynomial rings
whose second powers are not Ratliff-Rush closed:

\example
Let $I$ be the ideal in the polynomial ring $F[X,Y]$ over a field $F$
generated by
$$
Y^{22}, X^4Y^{18}, X^7Y^{15}, X^8Y^{14}, X^{11}Y^{11},
X^{14}Y^8, X^{15}Y^7, X^{18}Y^4, X^{22}.
$$
We will prove that this ideal is Ratliff-Rush closed,
but that $I^2$ is not.
Namely,
it is straightforward to verify that $X^{20} Y^{24}$ is not in $I^2$
but multiplies $I$ into $I^3$.

It remains to prove that $I$ is Ratliff-Rush closed.
It suffices to prove that every monomial in $\r I$ is in $I$.
As $I$ is primary to the maximal ideal $(X,Y)$,
it suffices to prove that every monomial in $\r I \cap (I : (X,Y))$ is in $I$.
But the only monomials in $I : (X,Y)$ not in $I$ are
$$
X^3Y^{21}, X^6Y^{17}, X^7Y^{14}, X^{10}Y^{13},
X^{13}Y^{10}, X^{14}Y^7, X^{17}Y^6, X^{21}Y^3.
$$
By symmetry it suffices to prove that the first four are not in $\r I$.
If $X^3 Y^{21}$ were in $\r I$,
then for some positive integer $n$,
$X^3 Y^{21} (Y^{22})^n \in I^{n+1}$.
But by the $X$-degree count then $X^3 Y^{21} (Y^{22})^n \in (Y^{22})^{n+1}$,
a contradiction.
If $X^6 Y^{17} \in \r I$,
then similarly $X^6 Y^{17} (Y^{22})^n \in (X^4 Y^{18}, Y^{22})^{n+1}$.
Thus $X^6 Y^{17} (Y^{22})^n \in (X^4 Y^{18})^a (Y^{22})^{n+1-a}$
for some non-negative integer $a \le n+1$.
By cancelling powers of $Y^{22}$,
then without loss of generality $X^6 Y^{17} (Y^{22})^n \in (X^4 Y^{18})^{n+1}$
for some $n$.
The $X$-exponents force $n = 0$, which is impossible.
Similarly,
$X^7 Y^{14}$ is not in $\r I$ for otherwise
$X^7 Y^{14} (Y^{22})^n \in (X^4 Y^{18})^a (X^7 Y^{15})^{n+1-a}$
for some non-negative integers $a$ and $n$ with $a \le n+1$.
By the $X$-degree count $n = 0$, which is impossible.
Finally,
$X^{10} Y^{13}$ is not in $\r I$ as otherwise
$X^{10} Y^{13} (Y^{22})^n \in (X^4 Y^{18}, X^7 Y^{15}, X^8 Y^{14})^{n+1}$
for some $n$.
By the $X$-degree count,
$n \le 1$,
and both cases are impossible.
This establishes that $I$ is Ratliff-Rush closed.

This last example also illustrates how difficult it is to prove
that an ideal is Ratliff-Rush closed.
Note also that in the example,
$I$ is not integrally closed.

Heinzer et al. constructed in [HJLS, (E3), page 386]
a Ratliff-Rush closed monomial ideal whose third power is not 
Ratliff-Rush closed.
Their example inspired the following construction
of monomial Ratliff-Rush closed ideals $I_n$
in the polynomial ring $F[X,Y]$
for which $(I_n)^n$ is not Ratliff-Rush closed:

\examples
Let $n \ge 3$ be an odd integer
and $I_n = (X^{3n-1}, X^{3n-4}Y^3, X^3 Y^{3n-4}, Y^{3n-1})$.
We will prove that $(I_n)^n$ is not Ratliff-Rush closed.
In fact,
we will prove that $(XY)^{(3n-1)n/2}$ is in $\r {(I_n)^n}$ but not in 
$(I_n)^n$.
To prove the inclusion
it suffices to prove (by Remark~\rmkrrred\ and symmetry)
that $(XY)^{(3n-1)n/2} (X^{3n-1})^{2m-1} \in (I_n)^{n + 2m - 1}$
for some positive integer $m$.
Let $m$ be such that $n = 2m+1$.
Then
$$
\eqalignno{
(XY)&^{(3n-1)n/2} (X^{3n-1})^{2m-1} =
(XY)^{6m^2 + 5m + 1} X^{12m^2 -2m-2}\cr &= (X^{18m^2+3m-1} Y^{9m+3}) 
\cdot (Y^{6m^2-4m-2}) = (X^{(6m-1)(3m+1)} Y^{3(3m+1)}) \cdot 
(Y^{(6m+2)(m-1)}) \cr &= (X^{3n-4}Y^3)^{3m+1} \cdot (Y^{3n-1})^{m-1} 
\in (I_n)^{4m} = (I_n)^{n + 2m -1}. \cr
}
$$
However,
$(XY)^{(3n-1)n/2}$ is not in $(I_n)^n$.
Otherwise there would exist non-negative integers $a, b, c$ and $d$
such that $a + b + c + d = n$ and
$$
(XY)^{(3n-1)n/2} \in (X^{3n-1})^a (X^{3n-4}Y^3)^b (X^3 Y^{3n-4})^c 
(Y^{3n-1})^d.
$$
The substitution $n = 2m+1$ yields the following equation for the 
$X$-exponents:
$$
(3m+1)(2m+1) = a(6m+2) + b(6m-1) + 3c
=2a(3m+1) +2b(3m+1) + 3(c-b),
$$
so that $(3m+1)(2m+1-2a-2b) = 3 (c-b)$.
Necessarily $2m+1-2a-2b$ is a multiple of $3$,
and hence $c-b$ is a multiple of $3m+1$,
which only holds if $c = b$ and hence $2m+1-2a-2b = 0$.
But the latter equation has no integer solutions.

\vskip3ex

 Another type of pathologies are ideals which are not Ratliff-Rush closed
but all of their powers are.
 For example the ideal $I=(X^4,X^3Y,XY^3,Y^4)$ in the polynomial ring $F[X,Y]$,
$F$ a field,
is not Ratliff-Rush closed
since $\r{I} =(X,Y)^4$,
but $\r{I^n}=I^n = (X,Y)^{4n}$ for every $n \ge 2$.
Here all the powers $I^n$, for $n \ge 2$, are also integrally closed. 

\vskip3ex

Finally we consider examples of ideals $I$
all of whose powers are Ratliff-Rush closed.
Let $gr_I(R)= \oplus_{n \ge 0} I^n/I^{n+1}$
denote the associated graded ring of $I$.
As was already proved in [RR, (2.3.1)],
if all the powers of $I$ are Ratliff-Rush closed,
then the depth of $gr_I(R)$ is positive.
The converse also holds
(for a proof,
see for example [HLS, (1.2)]):

\remark
\label{\rmkdepth}
All the powers of $I$ are Ratliff-Rush closed if and only if
the depth of $gr_I(R)$ is positive.

More generally,
in terms of the local cohomology of $G=gr_I(R)$
with respect to $G^+ = \oplus_{n \ge 1} I^n/I^{n+1}$,
for all integers $n$,
$$
{H^0}_{G^+}(G)_n = \r {I^{n+1}} \cap I^n /I^{n+1}.
$$
Thus ${I^{n+1}}$ is Ratliff-Rush closed if and only if
${H^0}_{G^+}(G)_n=0$.
Also, $I$ is Ratliff-Rush closed if and only if
there is no nonzero element in $G$
of degree zero which annihilates a power of $G^+$.

Thus the properties of the Ratliff-Rush closure of an ideal are a good
tool for getting the information on the depth of $G$
(see [S], [HM], [GR], [RV1],...),
and conversely,
when the depth of $G$ is positive,
we get information that all the powers of $I$ are Ratliff-Rush closed.

What are some ideals all of whose powers are Ratliff-Rush closed?
For example,
this holds if $I$ is an ideal
generated by a regular sequence in a Cohen-Macaulay ring.

Thus {\it parameter ideals in Cohen-Macaulay rings are all 
Ratliff-Rush closed.}
The Cohen-Macaulay assumption is necessary for this:

\example
\label{\exnotCM}
(Due to K.\ N.\ Raghavan, in [HJLS, Example 1.2])
Let $R$ be the subring $F[X,Y^2,Y^7,$
$X^2Y^5,X^3Y]$ of the polynomial ring $F[X,Y]$.
Then $I=(X,Y^2)R$ is a parameter ideal
(primary for the maximal ideal of $R$ of height two),
and $X^2Y^5 \in (I^2: I) \subseteq I$,
so $I$ is not Ratliff-Rush closed.
(As $Y^2 \cdot X^2 Y^5 = X^2 \cdot Y^7 \in XR$
and $X^2 Y^5 \not \in XR$,
then $R$ is not Cohen-Macaulay,
and $X, Y^2$ is not a regular sequence.)

Furthemore,
replacing an ideal generated by a regular sequence in a Cohen-Macaulay ring
by an {\it almost complete intersection
also fails to produce a Ratliff-Rush closed ideal.}
Namely,
Heinzer et al. in [HJLS, Theorem 2.2 and Corollary 2.3],
proved that if $I$ is a nonzero ideal of a Cohen-Macaulay local domain
which is an almost complete intersection ideal,
it is not necessarily Ratliff-Rush closed,
but it is minimal in its class,
i.e., for any ideal $J$ properly contained in $I$, $\r J \neq \r I$.
Here is an example of an almost complete intersection ideal
in a polynomial ring which is not Ratliff-Rush closed:

\example
\label{\exaci}
Let $I = (XY^5,X^6-Y^6,X^4Y^2-X^2Y^4)$ in the polynomial ring $F[X,Y]$.
Then $I$ is a zero-dimensional $3$-generated ideal,
so an almost complete intersection ideal.
 Note that $X^3Y^4 \not \in I$,
but that $X^3Y^4 \in I^3 : I^2$,
so that $X^3 Y^4 \in \r I$.

However,
zero-dimensional almost complete intersection monomial ideals
do satisfy the property that all of their powers are Ratliff-Rush closed:

\prop
\label{\propmuplusone}
Let $R$ be the polynomial ring
$F[X_1, \ldots, X_d]$ in $d$ variables over a field $F$.
Let $I$ be a $(d+1)$-generated zero-dimensional monomial ideal of $R$.
Then the associated graded ring $G$ of $I$ has positive depth.
In particular,
all the powers of $I$ are Ratliff-Rush closed.
\endb

\proof
Necessarily $d > 1$.
Note that $I = (X_1^{\alpha_1}, \ldots, X_d^{\alpha_d}, \gelt)$,
where for each $i$, $\alpha_i > \nu_i$.
Furthermore,
at least two of the $\nu_i$ are non-zero.
Suppose that $G$ has depth zero.
Then there exists an element $\rrelt \in I^n \setminus I^{n+1}$
(for some non-negative integer $n$)
such that $\rrelt (X_1, \ldots, X_d) \subseteq I^{n+1}$
and such that $\rrelt I \subseteq I^{n+2}$.

Let $I_0 = (X_1^{\alpha_1}, \ldots, X_d^{\alpha_d})$.
We will prove that whenever
$\rrelt I_0 \subseteq I^{n+2}$, $\rrelt \not \in I^{n+1}$ and $n \ge 0$,
then $\rrelt = \gelt \rselt$,
where $\rselt I_0 \subseteq I^{n+1}$.
Namely, by assumption
$$
\rrelt X_i^{\alpha_i} \in I^{n+2} \cap \left(X_i^{\alpha_i}\right)
\subseteq X_i^{\alpha_i} I^{n+1} + J^{n+2} \cap \left(X_i^{\alpha_i}\right),
$$
where $J = (X_1^{\alpha_1}, \ldots,
\widehat{X_i^{\alpha_i}}, \ldots, X_d^{\alpha_d}, \gelt)$.
As $\rrelt \not \in I^{n+1}$, necessarily
$$
\rrelt X_i^{\alpha_i} \in J^{n+2} \cap \left(X_i^{\alpha_i}\right)
\subseteq J^{n+2}.
$$
Write
$\rrelt X_i^{\alpha_i} =
\underline X^{\underline{\beta}} \left(X_1^{\alpha_1}\right)^{e_1} \cdots
\left(X_d^{\alpha_d}\right)^{e_d} \left(\gelt\right)^e$,
where the $e_j$ and $e$ are non-negative integers adding to $n+2$, $e_i = 0$,
and $\underline X^{\underline{\beta}}$ is a possibly unit monomial.
If $e = 0$,
then
$$
\rrelt \in
\left(X_1^{\alpha_1}, \ldots,
\widehat{X_i^{\alpha_i}}, \ldots, X_d^{\alpha_d}\right)^{n+2} : X_i^{\alpha_i}
=
\left(X_1^{\alpha_1}, \ldots,
\widehat{X_i^{\alpha_i}}, \ldots, X_d^{\alpha_d}\right)^{n+2},
$$
contradicting the hypotheses.
So necessarily $e > 0$.

\noindent Then from the equation
$\rrelt X_i^{\alpha_i} =
\underline X^{\underline{\beta}} \left(X_1^{\alpha_1}\right)^{e_1} \cdots
\left(X_d^{\alpha_d}\right)^{e_d} \left(\gelt\right)^e$,
necessarily for all $j \not = i$,
$\mu_j \ge \nu_j$.
As this holds for all $i$,
$\mu \ge \nu$ componentwise.
In particular,
for all $i$,
$\underline X^{\underline \mu-\underline\nu} X_i^{\alpha_i} =$
$\left(X_1^{\alpha_1}\right)^{e_1} \cdots $
$\left(X_d^{\alpha_d}\right)^{e_d} \left(\gelt\right)^{e-1} \in I^{n+1}$,
which proves the claim.

But this brings forth a contradiction:
by repeating this process,
$\rrelt$ is a product of $n+1$ copies of $\gelt$
times an element in $I : I_0 = R$,
so that $\rrelt \in I^{n+1}$.
\qed

This proposition cannot be extended to zero-dimensional monomial ideals
generated by $2$ more elements than the dimension of the polynomial ring:

\example
Let $I$ be the ideal $(X^{10},Y^5,XY^4,X^8Y)$
in the polynomial ring $F[X,Y]$.
Then $I$ is not Ratliff-Rush closed as
as $X^7Y^3 \in I^2 : I$ and
$X^7Y^3 \not \in I$.

The previous proposition,
together with Proposition~\propmuplusone,
shows in particular that the associated graded ring of
$I = (X^l, XY^{l-1}, Y^l)$ in $R = F[X,Y]$
has positive depth.
This answers [HJLS, (6.3), page 389].
In the next section we will need a more precise number:

\prop
\label{\propxyl}
Let $F$ be a field, $X, Y$ variables over $F$ and $R = F[X,Y]$.
Let $l$ be an integer $3$ or greater,
and let $I$ be the ideal $(X^l, XY^{l-1}, Y^l)$ in $R$.
Then the depth of the associated graded ring of $I$ is exactly $1$.
\endb

\proof
Let $G$ be the associated graded ring of $I$.
The image of $X^l$ in $G$ is a non-zero divisor as for all $n \ge 1$,
$$
\eqalignno{
I^{n+1} : X^l
&= \bigl(X^l I^n + (XY^{l-1},Y^l)^{n+1}\bigr) : X^l \cr
&= I^n + (Y^{l-1})^{n+1} (X,Y)^{n+1} : X^l
= I^n + (Y^{l-1})^{n+1}
(X,Y)^{(n+1){\rlap{\raise.7ex\hbox{.}}\hbox{-}} l}. \cr
}
$$
If $n+1 \le l$,
then the last ideal in the display is $I^n$ as $(l-1)(n+1) \ge ln$.
If instead $n+1 \ge l$,
then again the last ideal in the display is $I^n$
as
$(Y^{l-1})^{n+1} (X,Y)^{(n+1)- l} =
(Y^l)^{l-1} (XY^{l-1},Y^l)^{(n+1)- l} \subseteq I^{l-1+n+1-l} = I^n$.
Thus the depth of $G$ is at least $1$,
with $X^l$ being a non-zero divisor.
However,
$$
\eqalignno{
X \cdot (X^{l-1}Y) &\in I, \cr
Y^{l-2} \cdot (X^{l-1}Y) &\in I, \cr
X^l \cdot (X^{l-1}Y) &\in I^2 + (X^l), \cr
XY^{l-1} \cdot (X^{l-1}Y) &\in I^2 + (X^l), \cr
(Y^l)^{l-2} \cdot (X^{l-1}Y) &= (XY^{l-1})^{l-1} \in I^{l-1} + (X^l), \cr
}
$$
so that the non-zero image of $X^{l-1}Y$ in $G$
multiplies a power of the maximal ideal of $G$ into $X^lG$.
Thus the depth of $G$ is exactly $1$.
\qed

Another class of ideals whose associated graded ring has depth exactly 1 is:

\prop
\label{\propxyaal}
(Due to Shah-Swanson)
Let $d \ge 2$ be a positive integer
and let $R$ be the polynomial ring in variables $X_1, \ldots, X_d$
over a field $F$.
Let $l \ge 3$ be an integer and let
$I$ be the monomial ideal generated by all the monomials of degree $l$
except for the monomial $X_1^{l-1} X_2$.
Then the associated graded ring of $I$ has depth exactly $1$.
\endb

\proof
Note that for all $i > 1$,
$X_i (X_1^{l-1}X_2) = X_1 (X_1^{l-2}X_2X_i)$ is an element of $I$.
Thus
$$
\eqalignno{
(X_1, \ldots, X_d) (X_1^{l-1}X_2) & \subseteq I, \cr
I (X_1^{l-1}X_2) & \subseteq I^2 + (X_1^l). \cr
}
$$
Hence if $X_1^l$ is a non-zero divisor on the associated graded ring $G$,
then $G$ has depth exactly $1$.
So it remains to prove that $X_1^l$ is a non-zero divisor on $G$.
For this,
let $m$ be a monomial in $I^n \setminus I^{n+1}$
such that $m X_1^l \in I^{n+2}$.
We may write $m = m_0 m_1 m_2 \cdots m_n$,
where all the $m_i$ are monomials
with $m_1, \ldots, m_n$ some $l$-degree generators of $I$,
with possible repetitions.
By degree count $m_0$ has degree at least $l$,
and as by assumption $m \not \in I^{n+1}$,
necessarily $m_0$ is $X_1^{l-1}X_2$.
The assumption that $m X_1^l \in I^{n+2}$
implies that at least one of the $m_1, \ldots, m_n$,
say $m_1$,
is different from $X_1^l$.
But then by regrouping the variables,
$m_0 m_1 \in I^2$,
so that $m \in I^{n+1}$, contradicting the assumption.
\qed

\vfill\eject
\section{Ratliff-Rush closure and superficial elements}

In the study of the integral closure of ideals
(and of integral closures of powers of an ideal),
the notion of superficial elements is very important.
It enables one to drop dimension, for example.
In this section we review the elementary properties of superficial elements,
and show that in general there is no good notion of superficial elements
for the Ratliff-Rush closure.
At the end we look at special cases
where Ratliff-Rush closure behaves well
after reducing modulo a superficial (sufficiently general) element.

Recall the definition:
an element $a$ of an ideal $I$ is said to be {\it superficial} for $I$
if there exists a positive integer $c$ such that for all $n \ge c$,
$(I^n : a) \cap I^c = I^{n-1}$.

In particular,
if $a$ is a non-zero divisor in $R$,
$a\in I$ is superficial for $I$ if and only if
there exists a positive integer $c$ such that for all $n \ge c$,
$I^n : a = I^{n-1}$.
Notice that if $a$ is superficial for a regular ideal $I$, then $a 
\in I$, but $a
\not \in I^2$.

Integral closure generally behaves well after going
modulo a superficial element.
For example:

\prop
(Itoh [I, page 648])
If $I$ is an integrally closed $m-$primary ideal
in a Cohen-Macaulay local ring $(R,m)$ of dimension $d \ge 2$,
then at least after passing to a faithfully flat extension
there exists a superficial element $a \in I$
such that $IR/(a)$ is integrally closed in $R/(a)$.
\endb

In particular it follows that if $I$ is an integrally closed
$m-$primary ideal of a Cohen-Macaulay local ring $(R,m)$ of dimension
$d \ge 2$,
then $\r{I}=I$ and there exists a superficial element $a \in I$ such that
$IR/(a)$ is Ratliff-Rush closed in $R/(a)$.

However, superficial elements do not behave well for Ratliff-Rush 
closed ideals:
there exist many ideals $R$ all of whose powers are Ratliff-Rush closed,
yet for every superficial element $a \in I$,
$I/(a)$ is not Ratliff-Rush closed.

Before we give examples of this,
we prove a weaker version:

\prop
\label{\propdepthone}
There exist many ideals $R$ for which all the powers are Ratliff-Rush closed,
yet for every superficial element $a \in I$
there exists a positive integer $n$
such that $(I/(a))^n$ is not Ratliff-Rush closed in the ring $R/(a)$.
\endb

\proof
Let $I$ be any ideal for which the associated graded ring $G$
has depth exactly $1$.
By the result of Huckaba and Marley [HM],
every superficial element $a$ of $I$ is a non-zero divisor in $G$.
But $G$ modulo the image of $a$
is the associated graded ring for the ideal $I/(a)$ in the ring $R/(a)$,
and that has depth $0$.
Thus by Remark~\rmkdepth\ %
some power of $I/(a)$ is not Ratliff-Rush closed.
\qed

Such examples were provided in Propositions~\propxyl\ and~\propxyaal.
But we prove the stronger result for these two classes of examples:
it is actually $I/(a)$ itself
which is not Ratliff-Rush closed for all possible $a$.
First we prove this for examples from Proposition~\propxyl:

\prop
Let $R$ be the power series ring in variables $X$ and $Y$ over a field $F$.
Let $l > 2$ be an integer and let $I$ be the ideal $(X^l, XY^{l-1}, Y^l)$.
Then for every superficial element $a \in I$,
$I/(a)$ is not Ratliff-Rush closed.
\endb

\proof
Let $a = \alpha X^l + \beta XY^{l-1} + \gamma Y^l$,
with $\alpha, \beta, \gamma \in R$ and not all zero.

First note that $Y^l$ is not a superficial element as for every $n \ge 1$,
$X^{2+ln-l}Y^{l-2} \in I^{n-1} \setminus I^n$,
yet
$(X^{2+ln-l}Y^{l-2}) Y^l = (X^l)^{n-1} (XY^{l-1})^2 \in I^{n+1}$.
Also, $XY^{l-1}$ is not a superficial element as
$X^{ln-1}Y \in I^{n-1} \setminus I^n$,
but $(X^{ln-1}Y) (XY^{l-1}) = X^{ln} Y^l \in I^{n+1}$.

Suppose that $\alpha = 0$.
We just proved that then both $\beta$ and $\gamma$ are non-zero.
Set $b = X \sum_{i=0}^{l-1} (-1)^{l-1-i} \beta^i \gamma^{l-1-i} X^i Y^{l-i-1}$.
Then $ab = XY^{l-1} (\beta^l X^l + (-1)^{l-1} \gamma^l Y^l)$
is an element of $I^2$.
Thus for all $n \ge 1$,
$abX^{l(n-1)} \in I^{n+1}$,
yet
$bX^{l(n-1)} \in I^{n-1}\setminus I^n$.
Thus $a$ is not a superficial element and so necessarily $\alpha \not = 0$.
By localization at the maximal ideal $(X,Y)$
we can even conclude that $\alpha \not \in (X,Y)R$.

Now observe that $X^{l-1}Y \cdot XY^{l-1} = X^l Y^l \in I^2$,
$X^{l-1}Y \cdot Y^{l(l-1)} = (XY^{l-1})^{l-1} Y^l \in I^l$,
and so
$\alpha^{l-1} X^{l-1}Y \cdot X^{l(l-1)} \in
X^{l-1}Y (\beta XY^{l-1} + \gamma Y^l)^{l-1} + (a) \subseteq I^l + (a)$,
which proves that
$\alpha^{l-1}X^{l-1}Y \cdot I^{2l-3} \subseteq (a) + I^{2l-2}$.
Thus as $\alpha \not \in (X,Y)$,
it follows that $\alpha^{l-1} X^{l-1}Y \not \in I$,
so that $I/(a)$ is not Ratliff-Rush closed.
\qed

Similarly,
for examples from Proposition~\propxyaal,
$I/(a)$ itself is not Ratliff-Rush closed:

\prop
Let $d \ge 2$ be a positive integer
and let $R$ be the polynomial ring in variables $X_1, \ldots, X_d$
over a field $F$.
Let $l \ge 3$ be an integer and let
$I$ be the monomial ideal generated by all the monomials of degree $l$
except for the monomial $X_1^{l-1} X_2$.

Then for every superficial element $a \in I$,
$I/(a)$ is not Ratliff-Rush closed.
\endb

\proof
We denote by $J$ the subideal of $I$
generated by all the monomials of degree $l$
except for the monomials $X_1^{l-1} X_2$ and $X_1^l$
and we denote by $K$ the subideal of $I$
generated by $(X_1^{l+1},X_1^lX_2)$.

If $a \in I$ we may write
$a = \alpha X_1^l +b + c$ with $\alpha$ in $F$, $b$ in $J$ and $c$ in $K$.
Note that $X_1^{l-1}X_2 \not \in I$,
but $X_1^{l-1}X_2 (J+K) \subseteq I^2$.
Thus if $a$ is a superficial element for $I$,
as by [HM] its image is a non-zerodivisor in $gr_I(R)$,
necessarily $\alpha$ is a non-zero element of $F$.

%

For every superficial element $a \in I$, we prove now that
$\left(I^2 +(a)\right) : I \neq I$
and hence that $I/(a)$ is not Ratliff-Rush closed.
In fact $X_1^{l-1}X_2 \not \in I$, but $X_1^{l-1}X_2 I \subseteq I^2+(a)$ since
$$
\eqalignno{
X_1^{l-1}X_2 J &\subseteq I^2 \cr
X_1^{l-1}X_2 K &\subseteq I^2 \cr
X_1^{l-1}X_2 X_1^l&= \alpha^{-1}[ X_1^{l-1}X_2 a - X_1^{l-1}X_2 b] \in I^2+(a).
& \eqed
\cr
}
$$

Note that the same proof as in Proposition~\propdepthone\ shows:

\remark
Whenever $I$ is an ideal whose associated graded ring has depth at least $2$,
then for every superficial element $a \in I$ and every positive integer $n$,
$I^n$ and $I^n R/(a)$ are Ratliff-Rush closed.

Some examples of ideals with this property are the ideals $I$ in a Cohen-Macaulay ring
generated by a regular sequence of length $d \ge 2$.
(Note that these ideals need not be integrally closed.)
 
This class of examples cannot be extended to almost complete intersections
as was already shown in Example~\exaci.
It also cannot be extended to generic complete intersection ideals:

\example
Let $I = X(X^4,X^3Y,XY^3,Y^4)$ in the polynomial ring $F[X,Y]$.
Then {\it $I$ is a generic complete intersection which is not 
Ratliff-Rush closed}
as $X^3Y^2 \in (I^2 : I) \setminus I$.

\section{Initial ideals, associated primes and the number of generators}

In this section we concentrate on results showing
that the Ratliff-Rush closure does not behave well
with respect to several properties,
namely with respect to taking the initial ideals,
taking the associated primes,
and the number of generators.

First we analyze the Ratliff-Rush closure operation on
Borel-fixed and initial ideals.
We show that Borel-fixed ideals, even generic initial ideals,
need not be Ratliff-Rush closed,
and that Ratliff-Rush closure does not commute with taking the 
leading term ideal.
Furthermore,
the leading term ideal of a Ratliff-Rush closed ideal need not be 
Ratliff-Rush closed.

We also compare the sets of associated primes and the number of generators
of an ideal and of its Ratliff-Rush closure.
We show that there are no inclusion relations for the associated primes,
and similarly
that the number of generators of the Ratliff-Rush closure of an ideal
may be larger than the number of generators of the ideal.
We explore other connections between the associated
primes of $\r I$ and of various powers of $I$.
A modification of some corresponding arguments produces an answer to
Question~(1.6) in [HJLS]:
the minimal number of generators of $I$ need not be smaller
than the minimal number of generators of $\r I$.


The definition of Borel-fixed ideals can be found for example in [E, page 352].

Borel-fixed zero-dimensional ideals need not be Ratliff-Rush closed.
For example,
the ideal $I = (X^6,X^5Y^4,X^4Y^{11},X^3Y^{13}, X^2Y^{17},XY^{23},Y^{29})$
in the polynomial ring $F[X,Y]$
is Borel-fixed under the order $Y > X$.
However,
$I$ is not Ratliff-Rush closed as $X^4Y^9 \in I^2 : I$ and $X^4Y^9 \not \in I$.
Note that $I$ is even a generic initial ideal,
yet it is not Ratliff-Rush closed.

Furthermore,
there are many Borel-fixed zero-dimensional ideals
(even lex-segment ideals)
which are Ratliff-Rush closed but are not integrally closed.
For example,
in the polynomial ring $F[X,Y]$,
the ideal $I = (X^4,X^3Y,X^2Y^4,XY^5,Y^7)$
is zero-dimensional and lex-segment under the order $Y > X$,
but $I$ is not integrally closed because $X^2Y^3 \not \in I$
yet $(X^2Y^3)^2 = (X^3Y)(XY^5)\in I^2$.
One can calculate that $gr_I(R)$ has depth exactly one
so that by Proposition~\propdepthone,
for every superficial element $a$ in $I$,
$IR/(a)$ has some power which is not Ratliff-Rush closed.

We remark that for the Borel-fixed ideal
$I=(X^5, X^4Y, X^3Y^4, X^2Y^5, XY^8, Y^9)$,
the associated graded ring is Cohen-Macaulay of dimension $2$
but $I$ is not integrally closed.
Here $I$ and $I/(a)$ are Ratliff-Rush closed
for every superficial element $a \in I$.

Next we consider the behavior of the Ratliff-Rush closure under the
operation of taking the initial ideal.
We fix a polynomial ring and a monomial ordering.
The monomial ideal generated by all the leading terms of an ideal $I$
is denoted $\lt I$ (in the literature sometimes also in$\, I$).
Certainly $\lt I \subseteq \lt \r I$.
However,
we will show that $\r {\lt I}$ is not contained in $\lt \r I$,
and also that the other inclusion fails.

\example
Here is an ideal $I$ such that $\lt I$ is Ratliff-Rush closed,
$I$ is not Ratliff-Rush closed,
and $\r{\lt I} \subsetneq \lt \r I$.
Let $I = (XY^5,X^6-Y^6,X^4Y^2-X^2Y^4)$ in the polynomial ring $F[X,Y]$.
In Example~\exaci\ we showed that $X^3Y^4 \in \r I$.
Under the (reverse) lexicographic ordering,
$\lt I = (XY^5,X^6,X^4Y^2,Y^8)$.
It is easy to see that $\lt I$ is not Borel-fixed,
but is Ratliff-Rush closed.
In fact,
a computation by Macaulay2 shows that the associated graded ring of
$\lt I$ has positive depth,
so that $\r{\lt I} = \lt I$.
Then
$$
\r{\lt I} = (XY^5,X^6,X^4Y^2,Y^8)
\subsetneq (XY^5,X^6,X^4Y^2,Y^8,X^3Y^4) \subseteq \lt \r I.
$$

\example
Here is an ideal $I$ such that $I$ is Ratliff-Rush closed
but $\lt I$ is not,
and furthermore that $\lt \r I \subsetneq \r{\lt I}$.
Let $I = (X^7Y - X^2Y^5, X^5Y^2, X^2Y^5-XY^6, Y^7)$.
We verified by Macaulay2 and Singular that the associated graded ring of $I$
has positive depth,
so that $I$ (and all of its powers) are Ratliff-Rush closed.
It is easy to see that $\lt I = (X^7Y, X^5Y^2, X^2Y^5, Y^7)$,
and that $X^4Y^4 \in (\lt I)^2 : (\lt I)$.
Thus $X^4Y^4 \in \r {\lt I}$ but not in $\lt I = \lt \r I$.

This example also shows that the leading term ideal of a Ratliff-Rush closed
ideal need not be Ratliff-Rush closed.

Next we switch our attention to the sets of associated primes.
   From the definition of the Ratliff-Rush closure of an ideal $I$
it is clear that every associated prime ideal of $\r I$
is also associated to all high powers of $I$.
However,
there are no inclusions in general,
as the examples below show.

\example
\label{\exjarrah}
$\Ass {R/\r I}$ need not be contained in $\Ass {R/I}$,
even when $R/I$ is Cohen-Macaulay.

In~[J], Jarrah found a class of Cohen-Macaulay ideals $I$ for which
$\Ass {R/\overline I}$ is not contained in $\Ass {R/I}$.
The following $I$ is a modification of Jarrah's examples:
let $R$ be the polynomial ring $F[X,Y,Z]$,
and let $n \ge 2$ be a positive integer.
Set $J$ to be the ideal generated by all $X^i Y^{2n-i}$ except $i = n$.
Let $I = J \cap (X^n,Z) \cap (Y^n,Z)$.
Thus the set of associated primes of $I$ is $\{(X,Y), (X,Z), (Y,Z)\}$.
Clearly $R/I$ is Cohen-Macaulay.
Explicitly the generators of $I$ are
$$
I = J \cap (X^nY^n,Z)
= (X^i Y^{2n-i} Z | i \not = n) +
(X^n Y^{n+1}, X^{n+1} Y^n).
$$
By degree count,
$X^nY^n$ is not in the Ratliff-Rush closure of $I$.
Clearly $(X,Y) X^n Y^n \subseteq I$.
Also,
$X^n Y^n Z \in I^2 : I$,
so that $(X,Y,Z) X^n Y^n \subseteq \r I$.
This all proves that
$\r I : X^n Y^n = (X,Y,Z)$,
so that the maximal ideal $(X,Y,Z)$ is associated to $\r I$
but not to $I$.

(The key to this example is that $I$ localized at some
minimial prime is not Ratliff-Rush closed,
and that the other minimal components do not ``interfere'' with that.)

\example
$\Ass {R/I}$ need not be contained in $\Ass {R/\r I}$,
even if $R/\r I$ is Cohen-Macaulay.


But there exist also examples of monomial ideals $I$ in polynomial rings.
Namely, here is an example due to Ho\c sten:
$I = (X^4,X^3Y,XY^3,Y^4,X^2Y^2Z)$.
Ho\c sten constructed this ideal as an example of a monomial ideal $I$
for which $\Ass{R/I} \not \subseteq \Ass{R/I^2}$.
Clearly $Z$ is a zero-divisor modulo $I$.
In fact, the set of associated primes of $I$ equals $\{(X,Y), (X,Y,Z)\}$.
Note that $X^2Y^2 \in I^2 : I$,
thus $\r I = (X,Y)^4$,
and so $Z$ is not a zero-divisor modulo $\r I$,
so that $(X,Y,Z)$ is not an associated prime ideal of $\r I$.

\example
If $P$ is associated to $I^n$ for all $n$ sufficiently large,
$P$ need not be associated to $\r I$.

For this,
let $I$ be any three-generated height two prime ideal in
the polynomial ring $F[X,Y,Z]$.
(For example,
let $I = (X^3-YZ, Y^2-XZ, Z^2-X^2Y)$,
the kernel of the natural map $F[X,Y,Z] \to F[t^3,t^4,t^5]$.)
Huneke proved in~[Hn] that $(X,Y,Z)$ is associated to $I^n$ for all $n \ge 2$.
However,
$(X,Y,Z)$ is not associated to $\r I = I$.

Heinzer et al. asked in [HJLS, Question (1.6) (Q1)]
whether the minimal number of generators of a regular ideal
is always less than or equal to the minimal number of generators
of its Ratliff-Rush closure.
This seems to be the case also for all the examples so far in this paper.
However,
Ho\c sten's example can readily be modified to a counterexample to 
this question:

\example
Let $F$ be a field, $n \ge 2$ an integer,
$X, Y, Z_1, \ldots, Z_n$ variables over $F$,
$R = F[X,Y,Z_1, \ldots, Z_n]$,
and $I = (X^4, X^3 Y, X Y^3, Y^4) + (X^2 Y^2) (Z_1, \ldots, Z_n)$.
Then $\r I = (X,Y)^4$,
{\it the minimal number of generators of $I$ is $4 + n$,
and the minimal number of generators of $\r I$ is $5$.}

Furthermore,
here is even a zero-dimensional ideal illustrating this behavior:

\example
Let $R = F[X,Y,U,V]$,
$$
I = (X^4, X^3Y, XY^3, Y^4, X^2Y^2U, X^2Y^2V, U^4, V^4, XYU^2, XYV^2).
$$
Then $I$ is generated by $10$ monomials.
Its Ratliff-Rush closure is
$$
J = (X^4, X^3Y, X^2Y^2, XY^3, Y^4, U^4, V^4, XYU^2, XYV^2),
$$
generated by $9$ elements.
It is easy to prove that $J$ is in $\r I$,
and it suffices to prove that $J$ is Ratliff-Rush closed.
For this it suffices to prove that every monomial in $\r I$ is in $J$,
and as $I$ is zero-dimensional,
it suffices to prove that every monomial in $\r I \cap (J : (X,Y,U,V))$
is in $J$.
The only monomials in $J : (X,Y,U,V)$ not in $J$ are
$XY^2UV, X^2YUV, X^3U^3V^3, Y^3U^3V^3$.
By symmetry it suffices to prove
that the first and the last elements are not in $\r I$.
If $XY^2UV$ were in $\r I$,
then for some positive integer $n$,
$XY^2UV (Y^4)^n \in I^{n+1}$,
and hence by the $X$, $U$, and $V$-degree counts,
$XY^2UV (Y^4)^n \in (XY^3, Y^4)^{n+1}$,
which is impossible by the total $(X,Y)$-degree count.
Similarly,
if $Y^3 U^3 V^3 (Y^4)^n \in I^{n+1}$,
then $Y^3 U^3 V^3 (Y^4)^n \in (Y^4)^{n+1}$,
which is also impossible.

\vfill\eject
\section{The Ratliff-Rush reduction number }

Let $I$ be a regular ideal of a local ring $(R, m)$
and let $J$ be a minimal reduction of $I$.
We denote by
$$r_J(I):=min\{ n \ : \ I^{n+1}=JI^n \}$$
the {\it{reduction number }} of $I$ with respect to $J$.
If $I$ has a {\it {principal reduction}},
that is there exists $x \in I$
such that $I^{n+1}=x I^n$ for some integer $n$,
then $r_J(I)$ does not depend on $J$
(see [Hc], page 504).
If this is the case we write $r(I)$ instead of $r_J(I)$.
Notice that if $xR$ is a principal reduction of a regular ideal $I$,
then $x$ is a regular element of $R$.
It is known that $r(I) \le f_0 -1$,
where $f_0$ denotes the multiplicity of the fiber cone $F(I)=\oplus_{n\ge0}
I^n/mI^n$ (see [DGH], Corollary 5.3).
If $I$ has a principal reduction,
$F(I)$ is a one-dimensional graded ring which is not necessarily Cohen-Macaulay even
in the case $G=gr_I(R)$ is Cohen-Macaulay.

The notion of minimal reduction can be also given for
filtrations and the extension is clear in the case of the
Ratliff-Rush filtration. Since $\r {I^n} = I^n$ for large $n$,
a minimal reduction $J$ of $I
$ is a minimal reduction with respect to the Ratliff-Rush filtration.

We denote
by
$$
\widetilde r_J (I):= min\{ n \ : \ \r {I^{m+1}}=J \r{I^m} \ \ for\ \ m \ge \ n \}
$$
and we call it the {\it {Ratliff-Rush reduction number }} of $I$
with respect to $J$.
 It is not clear whether $\r {I^{n+1}}=J \r{I^n}$ for some
integer $n$ implies that $\r {I^{m+1}}=J \r{I^m}$ for every $m \ge n$.
We remark in fact that $\r I \r {I^n}$ is not necessarily $\r{I^{n+1}}$.
However for a regular ideal having principal reduction,
the Ratliff-Rush reduction number coincides with the least
integer $n$ such that
$\r{I^{n+1}}=x\r{I^n}$ where $(x)$ is a reduction of $I$,
which is the same behaviour that the $I-$adic filtration has.
In fact if $\r {I^{n+1}}=x \r{I^n}$ for some integer $n$,
then $\r {I^{m+1}}\subseteq \r{I^m}\subseteq (x)$ for 
every $m \ge n$,
hence $\r {I^{m+1}} \subseteq (x) \cap \r {I^{m+1}}=x\r {I^{m}}$.
Because $x\r {I^{m}} \subseteq \r {I^{m+1}}$,
the equality holds.

 Let $$ \r G:= \oplus_{n\ge0} \r {I^{n}}/\r{I^{n+1}}$$
the associated
graded ring to the Ratliff-Rush filtration of $I. $
It has a natural structure of graded algebra which has positive depth,
but it is not a standard algebra because we do
not necessarily have $\r {G_{n+1}}= \r{G_1} \r{G_n}$.
The following result shows that,
for a regular ideal having {\it {principal reduction}}, the Ratliff-Rush 
reduction number coincides with the least integer $t$ such that
$ \widetilde G_t \simeq \widetilde G_{t+1}$.
In particular this implies that if $I$ has a principal reduction,
$\widetilde r_J (I)$ again does not depend on $J$
and we write $\r r(I)$ instead of $\widetilde r_J (I)$.

The proof is inspired by Theorem 2.1.\ in [RV2].
We remark that the main tool used in the case of an
$\mm-$primary ideal of a
one-dimensional local ring was the information given by the
Hilbert function of $G$ and $\r G$,
the problem here is more complicated
because we have to replace this information with
different numerical invariants.

\prop
\label{\thmr} Let $I$ be a regular ideal of $(R, \mm)$
having principal reduction and let $t\ge 0$.
The following conditions are equivalent:
\item{a)}
$\widetilde G_t \simeq \widetilde G_{t+1}$.
\item{b)}
$ I\widetilde{I^t}+\widetilde{I^{t+2}}\subseteq
x\widetilde{I^t}+\widetilde{I^{t+2}}$.
\item{c)}
$\widetilde G_t \simeq \widetilde G_{n}$ for every $n \ge t$.
\vskip 2mm
\item{d)}
$\widetilde{I^{n+1}}= x \r {I^n}$
for some superficial element $x \in I$ and for every $n\ge t$.
\item{e)}
$\widetilde{I^{t+1}}= x \r {I^t}$ for some superficial element $x \in I$.
\endb

\proof
Since we have $\widetilde{I^{n+1}}:x=\widetilde{I^n}$ 
for every $n\ge 0$ and $x$ superficial element,
the multiplication by $x$ gives an injective map $$0\to
\widetilde{G}_n \ \overset{x}{\longrightarrow} \ \widetilde{G}_{n+1}$$ whose
cokernel is
$K_n:= \widetilde{G}_{n+1}/x\widetilde{G}_n=\widetilde{I^{n+1}}/(x\widetilde{I^n}+\widetilde{I^{n+2}})$.

It is clear now that a) implies b). Let us
prove that b) implies c). If $n\ge s$,
then $I^n = \r{I^n}$ and $I^{n+1}= x^{n+1-s}I^s$.
Since $x$ is regular in $R$,
we get
$\widetilde G_s= I^s/xI^s \simeq \widetilde G_{n}=$ for every $n \ge s$.

Hence we may assume $t \le n < s$ and we prove $\widetilde G_s
\simeq \widetilde G_{n}$.
By assumption
$I\widetilde{I^t}\subseteq x\widetilde{I^t}+\widetilde{I^{t+2}}$ and 
we claim that $I^s=x^{s-t}\widetilde{I^t}$.
We have
$$I^s=I^{s-t-1}I^{t+1}\subseteq I^{s-t-1}I\widetilde{I^{t}}\subseteq
I^{s-t-1}(x\widetilde{I^t}+\widetilde{I^{t+2}})\subseteq
xI^{s-t-1}\widetilde{I^t}+\widetilde{I^{s+1}}=xI^{s-t-1}\widetilde{I^t}+I^{s+1}.$$ 
If $s=t+1$ we are done by Nakayama. Otherwise $s> t+1$,
and
$$xI^{s-t-1}\widetilde{I^t}+I^{s+1}=xI^{s-t-2}I\widetilde{I^t}+I^{s+1}\subseteq
xI^{s-t-2}(x\widetilde{I^t}+\widetilde{I^{t+2}})+I^{s+1}\subseteq$$
$$ \subseteq x^2I^{s-t-2}\widetilde{I^t}+I^{s+1}\subseteq \dots
\subseteq x^{s-t}\widetilde{I^t}+I^{s+1}.$$ The claim follows again
by Nakayama.
\vskip 2mm
In particular for every $t \le n < s$,
$I^s=x^{s-n}\widetilde{I^n}$ and
$I^{s+1}=x^{s-n}\widetilde{I^{n+1}}$,
which imply c).

In fact $I^s=x^{s-t}\widetilde{I^t}\subseteq
x^{s-t-1}\widetilde{I^{t+1}}\subseteq \dots
x^{s-n-1}\widetilde{I^{n+1}}\subseteq \widetilde{I^{s }}=I^{s }. $
The second equality follows by using $I^{s+1}=xI^s$.

Let us finally prove that c) implies d). We may assume $n <s$.
If we prove that for every $t \le n < s$ we have $\widetilde{I^{n+1}}
\subseteq
x\widetilde{I^n} + \mm \widetilde{I^{n+1}}$,
then the conclusion follows by Nakayama.

Let $m$ be an integer such that $t \le m <s$.
By assumption we get $K_m=0$,
so it follows
$$\widetilde{I^{m+1} } \subseteq
x\widetilde{I^m}+\widetilde{I^{m+2}}.$$
We also have $K_{n }=0$ for $n > m$,
hence we get
$\widetilde{I^{m+2}} \subseteq
x\widetilde{I^m} + \widetilde{I^{m+3}} \subseteq \dots \dots
x\widetilde{I^m} + \widetilde{I^{s+1}}$
and we recall that $\widetilde{I^{s+1}}=x \widetilde{I^{s }} \subseteq
x \widetilde{I^{m+1}} \subseteq
\mm \widetilde{I^{m+1}}$.

Now d) implies e) is trivial and from e) we get a) because $K_t=0$.
\qed

A goal of this section is now to investigate the mutual relations
among $\r r (I)$
and
$r (I)$.

We give here a complete description in the case of a regular
ideal having principal reduction. We extend some results already known
in the case when $I$ is an $\mm-$primary ideal of a one-dimensional local
ring (see [RV2]).
Moreover the following Proposition, part i) recovers Proposition 
2.2., [DGH].

\prop \label{\reductionno}
Let $I$ a regular ideal of a local ring $(R, \mm)$
having principal reduction and let $s(I)$ be the smallest integer such that
$\r {I^n} = I^n$ for $n \ge s(I)$.
Then
\item{i)} $s(I)\le r(I)$,
\item{ii)} $\widetilde r(I) \le r(I)$,
\item{iii)} $\widetilde r(I) = r(I)$ if and only if $I^{r(I)} \not \subseteq (x)$
for some minimal reduction $(x)$ of $I$.

\endb

\proof
For brevity we put $r:=r(I)$.
Then $I^{i+1}=x I^i$ for every $i \ge r$,
from which it follows easily that for every $p\ge r$ and $t\ge0$,
$I^{r +t}=x^tI^{r }$.

Let $j$ be an integer, $j\ge r$, and let $t$ be a positive integer such that
$\widetilde{I^j}=I^{j+t}:I^t$.
Then
$$\widetilde{I^j}=I^{j+t}:I^t\subseteq
I^{j+t}:x^t=x^{j+t-r}I^r:x^t=x^{j-r}I^r\subseteq I^j,$$ so
that for every $j \ge r$, $\r I^j=I^j$.
This proves i) and ii).
Because $\r {I^{n }}\cap (x) = x\r {I^{n-1}}$ for every positive integer,
for the last assertion we remark that $\r {I^{n}}=x \r {I^{n-1 }}$ 
for some $n$
if and only if $\r{I^{n }} \subseteq (x)$.
Now the conclusion follows from i) because $\r {I^r}=I^r$.
\qed

It is clear that if $G$ has positive depth, then $\r r (I)=r(I)$ and $s(I)=0$.
The following example shows that the inequality ii) can be strict.

\example
Let $R=F[[t^4,t^5,t^{11}]]$ and $I$ the maximal ideal.
We have $I ^4=t^4 I ^3$,
so $r(I)=3$.
Now $I$ is Ratliff-Rush closed, but $I^2 \neq \r{I^2} =(t^8,t^9,t^{10},t^{11})$.
Moreover $\r{I^{n+1}}=t^4\r{I^{n }}$ for every $n \ge 2$,
hence $\r r (I)=2 < r(I)$.

\vskip 3mm

It is natural to ask if Proposition \reductionno\ can be extended to higher
analitic spread $\ell(I)$.
The following example shows that it cannot be extended to any
regular ideal $I$ having $\ell(I)=2$.

\example
As in the Example ~\exnotCM\ we consider the subring 
$R=F[X,Y^2,Y^7,X^2Y^5,X^3Y]$ of the polynomial ring $F[X,Y]$.
In this case $I=(X,Y^2)R$ is a parameter ideal,
in particular $r(I)=0$.
But $I$ is not Ratliff-Rush closed,
hence $\r r (I)\ge 1 > r(I)=0$.
We remark that $X, Y^2$ is not a regular sequence in $R$.

 In [DGH, Question (2.3) ] the authors asked
whether some of the good properties of the Ratliff-Rush filtration
established for regular ideals having principal reductions,
also hold for ideals having minimal reductions generated by regular sequences.

In this last part of this section assume $\ell(I)=2$ and
$J=(x,y)$ a minimal reduction of
$I$ generated by a regular sequence.

In the analogy with the ideal having principal reduction, we can
ask whether $\r {I^{n}} =I^n$ for $n \ge r_J(I)$.
In the case of an $\mm$-primary ideal of a local Cohen-Macaulay ring of dimension $d=2$,
Huneke proved that
if $I^{m+1}=JI^m$ and $I^{m+1}:x =I^m$,
then $I^{n+1}:x =I^n$ for all $n \ge m$ (see [HJLS], Proposition 4.3).
In particular we have $\r {I^n} = I^n$ for all $n \ge m$.
We do not know if $I^{m+1}=JI^m$ is enough to imply $\r {I^m} = I^m$.

In any case, we can prove the following result:

\prop
\label{\reduction}
With the above assumption, we have
$$ \widetilde r_J (I) \le r_J (I).$$
\endb

\proof
Put for short $r=r_J(I)$ and we prove
$\r {I^{r+1}}=J \r{I^r}$.
For a large $N$ we have $\r{I^r}= I^{r+N}:
(x^N,y^N)$ and $\r{I^{r+1}}= I^{r+N+1}: (x^N,y^N)$,
in particular $\r{I^{r+1}} = J^{N+1}I^r : (x^N,y^N)$
since $I^{r+n}=J^nI^r$ for all $n \ge 1$.

We claim that
$J^{N+1}I^r: (x^N,y^N) \subseteq J( I^{r+N}: (x^N,y^N))$.
To prove this,
let $a$ be an element in $J^{N+1}I^r: (a^N,b^N)$.
As $\{x,y\}$ is a regular sequence,
it is then an easy computation to prove that
$$
a \in JI^r + y (I^{r+N}: x^N) \cap x(I^{r+N}: y^N) \subseteq J(I^{r+N}: (x^N,y^N)).
\eqno \eqed
$$
 
We remark that the above result can be applied to any $\mm$-primary 
ideal $I$ of a $2-$dimensional local Cohen-Macaulay ring.
It is natural to ask the following question:

\question
Let $I$ an $\mm$-primary ideal of a local Cohen-Macaulay ring.
Is it always true that $\widetilde r_J (I) \le r_J (I) ?$
\endb

\vskip 3em
{\bgroup
\leftline{\bf References}
\bigskip

\font\eightrm=cmr8 \def\rm{\fam0\eightrm}
\font\eightit=cmti8 \def\it{\fam\itfam\eightit}
\font\eightbf=cmbx8 \def\bf{\fam\bffam\eightbf}
\rm
\baselineskip=10pt
\parindent=3.6em

\catcode`\.=11

\item{[DGH]}
M. D'Anna, A. Guerrieri and W. Heinzer,
Invariants of ideals having principal reductions,
{\it Comm. Algebra}, {\bf 29} (2001), 889-906.

\item{[E]}
D. Eisenbud,
Commutative Algebra with a View toward Algebraic Geometry,
Springer-Verlag, 1994.

\item{[GR]}
A. Guerrieri and M. E. Rossi,
Hilbert coefficients of Hilbert filtrations,
{\it J. Algebra}, {\bf 199} (1998), 40-61.

\item{[HJLS]}
W. Heinzer, B. Johnston, D. Lantz and K. Shah,
Coefficient ideals in and blowups of a commutative Noetherian domain,
{\it J. Algebra}, {\bf 162} (1993), 355-391.

\item{[HLS]}
W. Heinzer, D. Lantz and K. Shah,
The Ratliff-Rush ideals in a Noetherian ring,
{\it Comm. Algebra}, {\bf 20} (1992), 591-622.

\item{[Hc]}
S. Huckaba,
Reduction numbers for ideals of anallytic spread one, 
{\it J. Algebra}, {\bf 108} (1987), 503-512.

\item{[HM]}
S. Huckaba and T. Marley,
Hilbert coefficients for Hilbert filtrations,
{\it J. Algebra}, {\bf 199} (1997), 64-76.

\item{[Hn]}
C. Huneke,
The primary components of and integral closures of ideals in 3-dimensional
regular local rings,
{\it Math. Ann.}, {\bf 275} (1986), 617-635.


\item{[I]}
S. Itoh,
Hilbert coefficients of integrally closed ideals,
{\it J. Algebra}, {\bf 176} (1995), 638-652.


\item{[J]}
A. S. Jarrah,
Integral closures of Cohen-Macaulay monomial ideals,
preprint, 2001.


\item{[RR]}
R. J. \ Ratliff and D. E. Rush,
Two notes on reductions of ideals,
{\it Indiana\ Univ. Math. J.}, {\bf 27} (1978), 929-934.

\item{[RV1]}
M. E. Rossi and G. Valla,
A conjecture of J. Sally,
{\it Comm. Algebra}, {\bf 24 } (13) (1996), 4249-4261. 

\item{[RV2]}
M. E. Rossi and G. Valla,
On the Hilbert function of the Ratliff-Rush filtration, preprint, 2002.

\bigskip
\halign{\hskip 4em # \hfil & \hskip 8em # \hfil \cr
Dipartimento di Matematica & Department of Mathematical Sciences \cr
Universit\'a di Genova & New Mexico State University \cr
Via Dodecaneso 35 & Las Cruces, NM 88003-8001, USA \cr
19146 Genova, ITALY & {\tt iswanson@nmsu.edu} \cr
{\tt rossim@dima.unige.it} \cr
}

\egroup}

\end